\documentclass[11pt]{amsart}

\usepackage{amsmath}
\usepackage{amssymb}
\usepackage{tikz-cd}
\usepackage{enumerate}
\usepackage{comment}
\usepackage{bbm}
\usepackage{hyperref}

\allowdisplaybreaks

\theoremstyle{definition}
	\newtheorem{definition}{Definition}[section]
	\numberwithin{definition}{section}
	\newtheorem{remark}[definition]{Remark}
	\newtheorem{example}[definition]{Example}
	
	\newtheorem{assumption}[definition]{Assumption}
\theoremstyle{plain}
	\newtheorem{lemma}[definition]{Lemma}
	\newtheorem{proposition}[definition]{Proposition}
	\newtheorem{theorem}{Theorem}
	\newtheorem{corollary}[definition]{Corollary}

\numberwithin{equation}{section}

\DeclareMathOperator{\tr}{tr}
\DeclareMathOperator{\Id}{Id}

\newcommand{\map}[3]{#1\colon#2\rightarrow#3}
\newcommand{\deq}{\mathrel{\mathop:}=}
\newcommand{\xlto}[1]{\stackrel{#1}{\longrightarrow}}

\newcommand{\R}{\mathbb{R}}
\newcommand{\C}{\mathbb{C}}
\newcommand{\X}{\mathcal{X}}

\renewcommand{\H}{\mathrm{H}}
\newcommand{\Hc}{\H_{\mathrm{c}}}
\newcommand{\Hb}{\H_{\mathrm{b}}}
\newcommand{\Hcb}{\H_{\mathrm{cb}}}
\newcommand{\SU}{\mathrm{SU}}
\newcommand{\HC}{\mathbb{H}_{\mathbb{C}}}
\newcommand{\ad}{\mathrm{ad}}
\newcommand{\rk}{\mathrm{rk}}
\renewcommand{\HC}{\mathbb{H}_{\C}}
\newcommand{\Vol}{\mathrm{Vol}}
\newcommand{\delbar}{\overline{\partial}}
\renewcommand{\L}{\mathfrak}

\newcommand{\Ad}{\mathrm{Ad}}

\newcommand{\Omtilde}{\widetilde{\Omega}}
\newcommand{\NAH}{\mathrm{NAH}}
\newcommand{\gl}{\mathfrak{gl}}
\renewcommand{\sl}{\mathfrak{sl}}
\newcommand{\V}{\mathbf{V}}
\newcommand{\GL}{\mathrm{GL}}
\newcommand{\SL}{\mathrm{SL}}

\newcommand{\om}{\omega}
\newcommand{\Om}{\Omega}
\newcommand{\sig}{\sigma}

\newcommand{\ka}{\kappa}
\newcommand{\Ga}{\Gamma}
\renewcommand{\th}{\theta}
\newcommand{\vphi}{\varphi}

\newcommand{\gfr}{\mathfrak{g}}
\newcommand{\hfr}{\mathfrak{h}}
\newcommand{\kfr}{\mathfrak{k}}
\newcommand{\pfr}{\mathfrak{p}}

\newcommand{\ufr}{\mathfrak{u}}

\begin{document}

\title[Bounded cohomology, Higgs bundles, and Milnor--Wood inequalities]{Bounded cohomology, Higgs bundles, \\ and Milnor--Wood inequalities}

\author[T. Hartnick]{Tobias Hartnick}
\address{Institut f\"ur Algebra und Geometrie, KIT, Englerstra{\ss}e 2, 76128 Karlsruhe, Germany}
\email{tobias.hartnick@kit.edu}

\author[A. Ott]{Andreas Ott}
\address{Mathematisches Institut, Ruprecht-Karls-Universit\"at Heidelberg, Mathematikon, Im Neu-enheimer Feld 205, 69120 Heidelberg, Germany}
\email{aott@mathi.uni-heidelberg.de}


\begin{abstract}
	We explain how the generalized Milnor--Wood inequality for reductive representations of a cocompact complex-hyperbolic lattice into a Hermitian Lie group translates, under the non-abelian Hodge correspondence, into various kinds of Milnor--Wood inequalities for Higgs bundles. This clarifies the relation between the representation theoretic generalized Milnor--Wood inequality and the various different versions of Milnor--Wood inequalities for Higgs bundles that are known in the literature.
\end{abstract}

\maketitle

\tableofcontents

\section{Introduction and main results}
\label{SectionIntroductionAndMainResult}

\subsection{Milnor--Wood inequalities}

Over the last two decades there has been a strong interest in the study of representations of surface groups into higher rank real Lie groups, in particular in their moduli spaces and deformations, which nowadays goes under the name of higher Teichm\"uller theory. More recently, these studies have been extended to also include representations of fundamental groups of certain higher-dimensional manifolds (in particular, K\"ahler manifolds). Via the non-abelian Hodge correspondence (see Section \ref{SubSectionConstructingHiggs} below), originally due to Simpson, Donaldson and Corlette \cite{SimpsonLocalSystems,DonaldsonHarmonic,Corlette}, certain classes of representations correspond to various types of Higgs (principal) bundles. It is a natural and important problem to understand how geometric properties of objects on either side of the non-abelian Hodge correspondence relate to each other, see e.g.\,\cite{OSWW} for a recent example.

The present article studies a special instance of this correspondence. Namely, we will be concerned with bounds for characteristic numbers of representations and bundles which generalize the classical Milnor--Wood inequality for representations of surface groups into $\mathrm{PSL}_2(\R)$ \cite{Milnor, Wood}. Such inequalities have been established independently on the representation side and on the Higgs bundle side of the non-abelian Hodge correspondence, and it is one of our main goals in this paper to explain how they translate into each other.

We will work in the following setting: $M$ will be a closed complex-hyperbolic manifold with fundamental group $\Gamma$ and $G$ will be a connected simple non-compact real Hermitian Lie group with finite center. We will then consider reductive representations $\rho: \Gamma \to G$ and their corresponding $G$-Higgs principal bundles over $M$ (in the sense of Definition \ref{GHiggsPB} below). The version of the non-abelian Hodge correspondence that will be relevant for us in this paper was developed by Bradlow, Garcia-Prada, Gothen, and Mundet \cite{BGM,BGG} for the case of surface group representations, generalizing the work of Donaldson \cite{DonaldsonHarmonic}, Hitchin \cite{HitchinSelfDuality,LieGroups}, Corlette \cite{Corlette}, Simpson \cite{SimpsonVariations,SimpsonLocalSystems,SimpsonModuliI,SimpsonModuliII}, and by Chaput, Koziarz and Maubon for Higgs bundles over higher dimensional base manifolds \cite{KoziarzMaubon,KoziarzMaubon2,ChaputMaubon}; see also \cite{BGR,GarciaPradaGothenMundet/NAHC}. For the rest of this introduction we will always assume this setting.
 
On the representation side, one can establish Milnor--Wood type inequalities by employing techniques from bounded cohomology developed by Burger and Iozzi \cite{BIAnnouncement,BI2,BurgerIozzi/BoundedDifferentialForms,BI3}, and by Burger, Iozzi and Wienhard \cite{BurgerIozziWienhardHermitian,Surface}; see also \cite{Pozzetti}. These inequalities are most naturally formulated in terms of the so-called \emph{Toledo invariant} $T(\rho)$ of a representation $\rho: \Gamma \to G$, which in the generality discussed in this article was first introduced by Burger and Iozzi \cite{BIAnnouncement}, extending a classical definition of Toledo \cite{ToledoOld}. It was pointed out by Burger and Iozzi in \cite{BurgerIozzi/BoundedDifferentialForms} that the computation of the Gromov norm of the K\"ahler class of Hermitian Lie groups (due to Domic and Toledo \cite{DomicToledo/GromovNorm} and Clerc and \O rsted \cite{ClercOrsted/GromovNorm}) together with methods from bounded cohomology implies bounds for the absolute value of the Toledo invariant $T(\rho)$. We refer to these bounds as \emph{representation-theoretic inequalities of Milnor--Wood type}. These inequalities apply to fundamental groups $\Gamma$ of compact complex-hyperbolic manifolds of \emph{arbitrary dimension} (in fact, even to a much wider class of lattices) and to \emph{arbitrary Hermitian targets} $G$. Also, the proof of these inequalities (as recalled in Section \ref{SecMilnorWood} below) is uniform in that it does not involve any case-by-case analysis. Thus, on the representation side of the non-abelian Hodge correspondence we have a clean universal inequality of Milnor--Wood type, which applies in great generality. 

The situation on the Higgs bundle side is more complicated. With a given Higgs principal bundle one can associate vector bundles in various different ways, each of which admits its own Milnor--Wood type inequality in terms of certain characteristic numbers. Correspondingly, several different formulations of Milnor--Wood inequalities for Higgs bundles are found in the literature, see \cite{KoziarzMaubon,BGG,BGM,BGR,KoziarzMaubon2,ChaputMaubon}.

As we will explain in the next subsection, there are essentially two different ways to formulate Milnor--Wood type inequalities for Higgs bundles, which we call \emph{Milnor--Wood type inequalities in isotropy form}, respectively \emph{Milnor--Wood type inequalities in representation form}. The main goal of this article is  to clarify in what sense both types of inequalities appear as shadows of the universal representation-theoretic Milnor--Wood type inequality under the non-abelian Hodge correspondence. This helps to organize the zoo of Milnor--Wood inequalities for Higgs bundles, and explains why they are equivalent to each other and to the representation-theoretic version. Beyond that, our methods may also be of interest in their own right as they build a bridge between invariants from bounded cohomology and gauge-theoretic invariants.

\subsection{Isotropy form vs. representation form}
We now explain how the two different kind of Milnor--Wood type inequalities for Higgs bundles arise. Our starting point is the observation that the Toledo invariant of a reductive representation can be expressed as a characteristic number of two different (classes of) vector bundles associated to the corresponding Higgs principal bundles.

To begin with the simplest case, assume first that $\rho: \Gamma \to {\rm SU}(p,q)$ is a reductive representation into $G := {\rm SU}(p,q)$. Then the associated Higgs principal bundle (in the sense of  Definition \ref{GHiggsPB}) is a certain $K$-principal bundle, where $K = {\rm S}({\rm U}(p) \times {\rm U}(q))$. An obvious choice of representation for $K$ is the \emph{standard representation} on $\C^p \oplus \C^q$. Via this representation, every $G$-Higgs principal bundle gives rise to a pair of vector bundles ${\bf V}_{\rm std}^+$ and ${\bf V}_{\rm std}^-$ of ranks $p$ and $q$ respectively. A less obvious, but more canonical choice of representation of $K$ is given by the \emph{complexified isotropy representation} of $K$ on the complexified tangent space to $G/K$ at the basepoint. This representation splits into the $\pm i$-eigenspaces of the invariant complex structure on $G/K$, and hence every $G$-Higgs principal bundle gives rise to another pair of vector-bundles ${\bf V}_{\rm iso}^+$ and ${\bf V}_{\rm iso}^-$ whose ranks both equal $\dim G/K = 2pq$. One can now express the Toledo invariant and universal Milnor--Wood type inequality either in terms of the bundles ${\bf V}_{\rm std}^+$ and ${\bf V}_{\rm std}^-$ or in terms of the bundles  ${\bf V}_{\rm iso}^+$ and ${\bf V}_{\rm iso}^-$.

If $G$ is an arbitrary Hermitian Lie group (satisfying our standing assumptions), then the complexified isotropy representation is still defined, and hence we can still associate vector bundles ${\bf V}_{\rm iso}^+$ and ${\bf V}_{\rm iso}^-$ to a $G$-Higgs principal bundle. Instead of the standard representation of $ {\rm S}({\rm U}(p) \times {\rm U}(q))$ we can choose any representation $\sigma: K \to {\rm GL}(\C^p \oplus \C^q)$ which is \emph{admissible} in the sense of Definition \ref{DefAdmissibleRep} below. We then obtain associated vector bundles ${\bf V}_{\sigma}^+$ and ${\bf V}_{\sigma}^-$, which generalize the bundles ${\bf V}_{\rm std}^+$ and ${\bf V}_{\rm std}^-$ above. Again, we can formulate the universal inequality of Milnor--Wood type either in terms of the bundles  ${\bf V}_{\rm iso}^+$ and ${\bf V}_{\rm iso}^-$ or in terms of the bundles ${\bf V}_{\sigma}^+$ and ${\bf V}_{\sigma}^-$. This leads to \emph{Milnor--Wood type inequalities in isotropy form}, respectively \emph{Milnor--Wood type inequalities in representation form}.

\subsection{Statement of the results}

To state our results we introduce the following notation. Let $K<G$ be a maximal compact subgroup of $G$ and let $\X := G/K$ be the associated symmetric space. We write $\gfr = \kfr \oplus \pfr$ for the corresponding Cartan decomposition, where $\pfr = T_{eK}\X$, and denote by ${\rm ad}_{\kfr}^{\pfr}: \kfr \to {\rm End}(\pfr)$ the infinitesimal isotropy representation of $\kfr$ on $\pfr$. Since $G$ is Hermitian, the space $\X$ admits a $G$-invariant complex structure $I_{\X}$, and we fix one of the two possible choices of $I_{\X}$ (equivalently, an orientation on $\X$) once and for all. The space $\Omega(\X)^G$ of $G$-invariant $2$-forms on $\X$ is one-dimensional, generated by the \emph{canonical $2$-form} $\om_{G}^{\mathrm{can}}$ which is uniquely determined by the fact that at the basepoint $o = eK$,
\[
(\om_{G}^{\mathrm{can}})_o(X,Y) = \tr_{{\rm End}(\pfr)} \left( (I_{\X})_o \circ \ad_{\kfr}^{\pfr}([X, Y]) \right) \quad  (X, Y \in \pfr = T_o\X).
\]
If $\omega \in \Omega(\X)^G$ is non-zero then exactly one of $\omega(\cdot, I_{\X}\cdot)$ and $-\omega(\cdot, I_{\X}\cdot)$ is a Riemannian metric on $\X$, and we denote the minimal holomorphic sectional curvature of this Riemannian metric by $K^{\min}_{\om}$.

Now denote by $\omega_M$ the K\"ahler form on $M$ which pulls back to the hyperbolic K\"ahler form on the universal cover of $M$. We denote by
\[
\Vol(M) \deq \int_{M} \om_{M}^{n}
\]
the hyperbolic volume of $M$. Given a vector bundle ${\bf V}$ over $M$, we define
\[
	{\bf deg}( \V) \deq = \int_{M} c_{1}( \V) \wedge \om_{M}^{n-1},
\]
where $c_1(\V) = c_1(\det(\V))$ denotes the first Chern class of $\V$.

The following theorem is established in Corollary \ref{MilnorWoodAdjoint}.

\begin{theorem}[Universal Milnor--Wood inequality, isotropy form] The degrees of the bundles ${\bf V}_{\rm iso}^+$ and ${\bf V}_{\rm iso}^-$ as defined above satisfy the inequality
\[
\bigl|{\bf deg}({\bf V}_{\rm iso}^+)  \bigr| = \bigl|{\bf deg}({\bf V}_{\rm iso}^-)  \bigr| \le \frac{{\rm Vol(M)} \cdot \rk(G)}{4\pi \cdot |K^{\rm min}_{\om_{G}^{\rm can}}|}.
\]
\end{theorem}
The constant $K^{\rm min}_{\om_{G}^{\rm can}}$ can be expressed in terms of root data of the Lie algebra of $G$. For example for $G = {\rm SU}(p,q)$ we have $K^{\rm min}_{\om_{G}^{\rm can}}= -\frac{1}{p+q}$, and thus we obtain
\[
\bigl|{\bf deg}({\bf V}_{\rm iso}^+)  \bigr| = \bigl|{\bf deg}({\bf V}_{\rm iso}^-)  \bigr| \le \frac{\rm Vol(M)}{4\pi} \cdot (p+q) \cdot \min\{p,q\}.
\]

To state the second version of the Milnor--Wood type inequality we observe that for every admissible representation $\sigma: G \to {\rm GL}(V)$ in the sense of Definition \ref{DefAdmissibleRep}, for all $X, Y \in \pfr$ the complex trace $\tr_\C(d\sig(I) \circ d\sig([X,Y]))$ is actually a real number and there exists a unique non-zero $G$-invariant $2$-form $\omega^\sigma$ on $\X$ such that
\[(\om^{\sig})_o(X, Y) = \tr_\C((I_\X)_o \circ d\sig([X,Y]))\]
for all $X, Y\in \pfr$ (see Lemma \ref{LemmaExtendTraceForm}).

The following theorem is established in Corollary \ref{MilnorWoodAdmissible}.

\begin{theorem}[Universal Milnor--Wood type inequality, representation form] 
The degrees of the bundles ${\bf V}_{\sigma}^+$ and ${\bf V}_{\sigma}^-$ associated with an admissible representation $\sigma$ satisfy the inequality
\[
\bigl| \deg(\V_{\sig}^{+}) \bigr| = \bigl| \deg(\V_{\sig}^{-}) \bigr| \le \frac{\Vol(M)  \cdot \rk(G)}{2\pi \cdot |K^{\rm min}_{\om^{\sig}}|}.
\]
\end{theorem}
Again, the constant $ |K^{\rm min}_{\om^{\sig}}|$ can be computed explicitly for a given admissible representation $\sigma$. For the standard representation $\sigma: {\rm SU}(p,q) \to {\rm GL}(\C^p \oplus \C^q)$ we obtain$
|K^{\rm min}_{\om^{\sig}}| = 2$, and hence
\begin{equation*}
	\bigl|{\bf deg}({\bf V}_{\rm std}^{+})  \bigr| = \bigl|{\bf deg}({\bf V}_{\rm std}^{-})  \bigr|  \le \frac{{\rm Vol(M)}}{4\pi} \cdot {\rm min}\{p,q\}.
\end{equation*}
In the special case where $M = \Sigma_g$ is a closed hyperbolic surface of genus $g \geq 2$ we have $ \frac{{\rm Vol(M)}}{4\pi} = g-1$, which is in accordance with the Milnor--Wood inequality obtained in \cite{BGG}.

\bigskip

\noindent \textbf{Acknowledgements.} We thank Kloster Mariastein for their hospitality and excellent working conditions. We would also like to thank Marc Burger, Oscar Garc\'ia-Prada, Alessandra Iozzi, Vincent Koziarz, Ignasi Mundet i Riera, Julien Maubon, Ana Pe\'on-Nieto, Andy Sanders, Roberto Rubio, Anna Wienhard, and Chris Woodward for discussions and suggestions, and for pointing out some mistakes.\\
	A.\,O.\, was supported by grant TH-01 06-1 from ETH Z\"urich, the European Research Council under ERC-Consolidator Grant 614733 ``Deformation Spaces of Geometric Structures'', and by the Priority Program 2026 ``Geometry at Infinity'' of the Deutsche Forschungsgemeinschaft (DFG, German Research Foundation) under DFG grant 340014991.  He further acknowledges funding by the Deutsche Forschungsgemeinschaft (DFG, German Research Foundation) -- 281869850 (RTG 2229). He was also supported by the Deutsche Forschungsgemeinschaft (DFG, German Research Foundation) under Germany's Excellence Strategy EXC 2181/1 - 390900948 (the Heidelberg STRUCTURES Excellence Cluster).

\section{Toledo invariant and Milnor--Wood type inequality}\label{SecMilnorWood}

\subsection{Continuous bounded cohomology and the Gromov norm}

Let $G$ be a locally compact group. The \emph{continuous cohomology} $\Hc^{\bullet}(G;\R)$ of $G$ with real coefficients \cite{BorelWallach} is the cohomology of the complex $(C_{\mathrm{c}}^{\bullet}(G;\R),\delta)$, where
\[
C_{\mathrm{c}}^{n}(G;\R) \deq C(G^{n+1},\R)^{G}
\]
denotes the space of continuous functions $G^{n+1}\to\R$ which are invariant under the diagonal $G$-action on $G^{n+1}$, and
\[
(\delta f) (g_{0},\ldots,g_{n}) \deq \sum_{i=0}^{n} (-1)^{i} \cdot f(g_{0},\ldots,\widehat{g_{i}},\ldots,g_{n}).
\]

Likewise, the \emph{continuous bounded cohomology} $\Hcb^{\bullet}(G;\R)$ of $G$ with real coefficients \cite{MonodBook} is defined to be the cohomology of the subcomplex $(C_{\mathrm{cb}}^{\bullet}(G;\R),\delta)$ of $(C_{\mathrm{c}}^{\bullet}(G;\R),\delta)$, where $C_{\mathrm{cb}}^{n}(G;\R) \subset C_{\mathrm{c}}^{n}(G;\R)$ denotes the space of continuous bounded $G$-invariant functions $G^{n+1}\to\R$.

Given a class $\alpha \in \Hcb^{n}(G;\R)$ in the bounded cohomology of $G$, its \emph{Gromov norm} is defined as
\[
\|\alpha\| \deq \inf\bigl\{ \|c\|_{\infty} \,\big|\, c \in \alpha \bigr\},
\]
where the infimum is taken over all cocycles $c$ representing the class $\alpha$, and $\|c\|_{\infty}$ denotes the sup-norm of the bounded function $\map{c}{G^{n+1}}{\R}$.

\subsection{Hermitian symmetric spaces and the bounded K\"ahler class}
\label{SubSectionSymmetricSpaces}

We recall basic facts about Hermitian symmetric spaces and introduce some notation. Let $G$ be a simple non-compact Lie group with finite center. Fix a maximal compact subgroup $K$ of $G$ and denote by $\gfr = \kfr \oplus \pfr$ the Cartan decomposition of the corresponding Lie algebras. We will write $K_{\C}$, $G_{\C}$, $\kfr_{\C}$, $\pfr_{\C}$ and $\gfr_{\C}$ for the complexifications of $K$, $G$, $\kfr$, $\pfr$ and $\gfr$ respectively.

We denote by $\X_G = G/K$ the associated symmetric space and fix an orientation for $\X_G$. We also choose the identity coset $o = eK$ as a basepoint in $\X_G$ and identify the tangent space $T_{o}\X_G$ of $\X_G$ at $o$ with $\pfr$. The Lie group $G$, and also the symmetric space $\X_G$, are called \emph{Hermitian} if $\X_G$ carries a $G$-invariant complex structure. In this case, there is a unique $G$-invariant complex structure $I_G$ which is compatible with the chosen orientation. From now on we will always consider $\X_G$ as a complex manifold with complex structure $I_G$. Note that there exists a unique element $I \in \kfr$ which is central in $\kfr$ and satisfies
\begin{equation} \label{DefI}
	(I_{G})_o(X) = \ad(I)(X) = [I,X]
\end{equation}
for all $X \in \pfr$.

The space $\Om^{2}(\X_{G})^{G}$ of $G$-invariant $2$-forms is one-dimensional. For any such form $\om \in \Om^{2}(\X_{G})^{G}$ we define
\begin{equation} \label{EquationCompatibleTriple}
	g_{\omega}(\,\cdot,\,\cdot) = \om(\,\cdot,I_{G}\,\cdot).
\end{equation}
If $\om$ is non-zero, then either $g_{\om}$ or $-g_{\om}$ is a Riemannian metric, and hence $\om$ or $-\om$ is a K\"ahler form. We then denote by $K^{\min}_{\om}$ the minimal holomorphic sectional curvature of the Riemannian metric $\pm g_{\om}$ (see e.g.~\cite{Ballmann}). Note that $\om$ is uniquely determined up to sign by the invariant $K^{\min}_{\om}$.

Recall that the comparison map $\Hcb^2(G; \R) \to \Hc^2(G; \R)$ is an isomorphism (see \cite{BurgerIozzi/BoundedDifferentialForms}, Section 5). Composing the van Est isomorphism $\Omega^2(\X_G)^G \to \Hc^2(G; \R)$ (see \cite{Guichardet/Cohomologie}, Ch.\,III \S7) with the inverse of the comparison map we obtain the \emph{bounded van Est isomorphism}
\begin{equation} \label{MapVanEst}
	\map{\iota}{\Om^{2}(\X_G)^{G}}{\Hcb^{2}(G;\R)}.
\end{equation}
Following Burger and Iozzi \cite{BurgerIozzi/BoundedDifferentialForms}, Section 5, we associate with every non-zero $G$-invariant $2$-form $\om \in \Om^{2}(\X_{G})^{G}$ the associated \emph{bounded K\"ahler class}
\[
\ka_{\om}^{b} \deq \iota(\om) \in \Hcb^{2}(G;\R).
\]
Note that $\ka_{\om}^{b}$ generates $\Hcb^{2}(G;\R)$. Fixing a reference point $\overline{x} \in \X_G$, we can represent the class $\ka_{\om}^{b}$ by the cocycle 
\[
\map{c_{\om}}{G^{3}}{\R}, \quad c_{\om}(g_{0},g_{1},g_{2}) \deq \int_{\bigtriangleup(g_{0}.\overline{x},g_{1}.\overline{x},g_{2}.\overline{x})} \om,
\]
where $\bigtriangleup(x_{0},x_{1},x_{2})$ denotes, for any three points $x_{0}, x_{1}, x_{2} \in \X_G$, the geodesic triangle in $\X_G$ with vertices $x_{0}, x_{1}, x_{2}$, taken with respect to some fixed $G$-invariant Riemannian reference metric on $\X_G$.

\subsection{Toledo invariant}

Let $G$ and $H$ be simple Hermitian Lie groups  (always assumed non-compact and with finite center) with associated Hermitian symmetric spaces $(\X_{G}, I_G)$ and $(\X_{H}, I_H)$. Let $\om_{G} \in \Om^{2}(\X_{G})^{G}$ be an arbitrary non-zero $G$-invariant  $2$-form on $\X_G$, and let $\omega_H$ be the unique $H$-invariant K\"ahler form on $\X_H$ normalized to $K_{\omega_H}^{\rm min} = -1$.

Fix a cocompact lattice $\Ga < H$, and consider a representation $\map{\rho}{\Ga}{G}$. The representation $\rho$ gives rise to a map
\begin{equation} \label{MapExtendedToledoMap}
	T_{b}(\rho)\colon \Om^{2}(\X_G)^{G} \xlto{\iota} \Hcb^{2}(G;\R) \xlto{\rho^{\ast}} \Hb^{2}(\Ga;\R) \xlto{T_{b}} \Hcb^{2}(H;\R)
\end{equation}
which is defined as the composition of the bounded van Est isomorphism \eqref{MapVanEst}, the natural pull-back map $\map{\rho^{\ast}}{\Hcb^{2}(G;\R)}{\Hb^{2}(\Ga;\R)}$, and the bounded transfer map $\map{T_{b}}{\Hb^{2}(\Ga;\R)}{\Hcb^{2}(H;\R)}$ (see \cite{BurgerIozzi/BoundedDifferentialForms}, Sections 2.1 and 2.2 for details). Since $\Hcb^{2}(H;\R) \cong \R \cdot \ka_{\om_{H}}^{b}$ the following definition makes sense.

\begin{definition}
	Fix a representation $\map{\rho}{\Ga}{G}$ and a non-zero invariant $2$-form $\om_{G} \in \Om^{2}(\X_{G})^{G}$. We define the \emph{Toledo invariant} $T_{\om_{G}}(\rho)$ of $\rho$ with respect to $\om_{G}$ by the relation
\begin{equation} \label{EquationDefinitionToledoInvariant}
	T_{b}(\rho)(\om_{G}) = T_{\om_{G}}(\rho) \cdot \ka_{\om_{H}}^{b}.
\end{equation}
\end{definition}

This definition is taken from Burger and Iozzi \cite{BurgerIozzi/BoundedDifferentialForms}, up to composition with the bounded van Est isomorphism and the fact that we do not assume any normalization condition for $\om_{G}$. Consequently, our invariant $T_{\om_{G}}(\rho)$ will depend on $\om_{G}$. We warn the reader, that different choices of normalization exist in the literature.

\subsection{Milnor--Wood inequality}

We keep the assumptions and notation from the previous section. Our knowledge of the Gromov norm of the bounded K\"ahler class due to Domic and Toledo \cite{DomicToledo/GromovNorm} and Clerc and {\O}rsted \cite{ClercOrsted/GromovNorm} immediately yields the following variant of the Milnor--Wood inequality for the Toledo invariant, which is due to Burger and Iozzi. 

\begin{lemma}[Milnor--Wood inequality (\cite{BurgerIozzi/BoundedDifferentialForms}, Thm.\,7)]\label{LemmaMW}
	Let $G$ and $H$ be simple Hermitian Lie groups with associated symmetric spaces $\X_{G}$ and $\X_{H}$ and let $\om_{G} \in \Om^{2}(\X_{G})^{G}$ be a non-zero invariant $2$-form on $\X_{G}$. Fix a cocompact lattice $\Ga < H$, and consider a representation $\map{\rho}{\Ga}{G}$. Then the Toledo invariant $T_{\om_{G}}(\rho)$ of $\rho$ with respect to $\om_{G}$ satisfies
\[
\bigl| T_{\om_{G}}(\rho) \bigr| \le \frac{1}{|K^{\rm min}_{\omega_G}|} \cdot \frac{\rk(G)}{\rk(H)}.
\]
\end{lemma}

\begin{proof}
	First of all we note that multiplying $\omega_{{G}}$ by a scalar $\lambda > 0$ has the effect of multiplying $K^{\min}_{\omega_G}$ by a factor of $1/\lambda$ (see \cite{ClercOrsted/GromovNorm}, p.\,274 for details). Since the Toledo invariant $T_{\om_{G}}(\rho)$ depends linearly on $\om_{G}$, it follows that we may without loss of generality assume that $\omega_G$ is K\"ahler and normalized to $K^{\min}_{\omega_G} = -1$.

With this normalization understood, by work of Domic and Toledo \cite{DomicToledo/GromovNorm} and Clerc and {\O}rsted \cite{ClercOrsted/GromovNorm} the Gromov norms of the bounded K\"ahler classes $\ka_{\om_{G}}^{b}$ and $\ka_{\om_{H}}^{b}$ are given by
\[
\| \ka_{\om_{G}}^{b} \| = \pi \cdot \rk(G), \quad \| \ka_{\om_{H}}^{b} \| = \pi \cdot \rk(H).
\]
Hence the defining relation \eqref{EquationDefinitionToledoInvariant} together with the fact that the map $T_{b} \circ \rho^{\ast}$ in \eqref{MapExtendedToledoMap} does not increase the Gromov norm yields
\begin{eqnarray*}
|T_{\om_{G}}(\rho)| \cdot \pi \cdot \rk(H) \; &=& \;  |T_{\om_{G}}(\rho)| \cdot \| \ka_{\om_{H}}^{b} \| \; = \; \| T_{\om_{G}}(\rho) \cdot \ka_{\om_{H}}^{b} \|  \; = \;  \| T_{b}(\rho)(\om_{G}) \|\\
&=& \;  \| (T_{b} \circ \rho^{\ast} \circ \iota)(\om_{G}) \| \; \leq \; \| \iota(\om_{G}) \|\; = \; \| \ka_{\om_{G}}^{b} \| \\ &=& \;\pi \cdot \rk(G).
\end{eqnarray*}
Dividing both sides by $\pi \cdot \rk(H)$ gives the desired inequality.
\end{proof}

For later use let us also point out how the constant $K^{\rm min}_{\omega_G}$ can be computed in practice. We first recall (see e.g. \cite{BurgerIozziWienhardHermitian}) that if $\L a\subset \pfr$ is a maximal abelian subalgebra, then the root system of $\gfr$ with respect to $\L a$ is of type $C_n$ or $BC_n$ according to whether $\X_G$ is of tube type or not. Thus there exist a basis $(\xi_1, \dots, \xi_r)$ of $\L a^*$ and integers $a,b$ such that the roots are given by $\pm 2\xi_j$, $\pm \xi_j \pm \xi_k$ ($j \neq k$) and possibly $\xi_j$ with respective multiplicities $1$, $a$ and $b$. The integers $a$, $b$ and $r$ determine $\gfr$ uniquely. A fundamental invariant of $\gfr$ is the number
\begin{equation}\label{pg}
p_\gfr := (r-1)\cdot a + b +2.
\end{equation}
The following is an immediate consequence of \cite{ClercOrsted/GromovNorm}:
\begin{lemma} Let $\om_{G} \in \Om^{2}(\X_{G})^{G}$ be a non-zero $G$-invariant $2$-form on $\X_G$ and $I \in \kfr$ as in \eqref{DefI}. Then for every non-zero $X \in T_o\X_G =\pfr$ we have
\begin{equation}\label{KMinFormula}
|K^{\rm min}_{\omega_G}| = \frac{1}{p_\gfr} \cdot \frac{|\tr({\rm ad}(X)^2)|}{|(\omega_G)_o(X,[I, X])|}.
\end{equation}
\end{lemma}
\begin{proof} Replacing $\omega_G$ by $-\omega_G$ if necessary, we may without loss of generality assume that $\omega_G$ is a $G$-invariant K\"ahler form, i.e., there exists  a $G$-invariant Hermitian metric $h = g + i\, \omega_G$ on $\X_G$. We can realize $\X_G$ as a bounded symmetric domain and denote by $h_{\rm Berg} = g_{\rm Berg} + i\, \omega_{Berg}$ the corresponding Bergman metric. Since the space of invariant metrics on $\X_G$ is one-dimensional we have $h = \lambda \cdot h_{\rm Berg}$ for some $\lambda \in \R$. To compute $\lambda$ we recall two basic facts from \cite{ClercOrsted/GromovNorm}: Firstly, $h_{\rm Berg}$ has minimal holomorphic sectional curvature $- \frac 2{p_{\gfr}}$. Secondly, the restriction of $g_{\rm Berg}$ to $\pfr$ is given by $\frac 1 2 \kappa$, where $\kappa \in \bigwedge^2 \pfr^*$  denotes the restriction of the Killing form $\kappa(X,Y) = {\rm tr}({\rm ad}(X){\rm ad}(Y))$ to $\pfr$.  Now for every $X \in \pfr$ we have
\[
\lambda =  \frac{h_o(X, X)}{(h_{\rm Berg})_o(X,X)} = \frac{g_o(X, X)}{(g_{\rm Berg})_o(X,X)} = \frac{g_o(X,X)}{\frac 1 2 \tr({\rm ad}(X)^2)},
\]
whence
\[
K^{\rm min}_{\omega_G}= - \frac{1}{\lambda} \cdot \frac 2{p_{\gfr}} = - \frac{\frac 1 2 \tr({\rm ad}(X)^2)}{g_o(X,X)} \cdot  \frac 2{p_{\gfr}} = -\frac{1}{p_{\gfr}} \cdot \frac{\tr({\rm ad}(X)^2)}{(\omega_G)_o(X,[I,X])}.
\]
\end{proof}

\section{The non-abelian Hodge correspondence revisited}
\label{SectionHodgeCorrespondence}

\subsection{Notation and preliminaries}
\label{SecToledoHiggsPrelims}

Keeping the notation from the previous section, we now specialize to the case where $H = {\rm SU}(1, n)$. In this case, the associated symmetric space $(\X_{H},I_{H})$ can be identified with complex hyperbolic $n$-space $\HC^{n}$ together with the standard complex structure. We denote by $\om_{\HC^{n}}$ the unique invariant K\"ahler form on $\HC^{n}$ normalized to minimal holomorphic sectional curvature $K_{\om_{\HC^{n}}}^{\min} = -1$.

By a compact complex-hyperbolic manifold we understand a manifold $M$ of the form $M = \Ga \backslash \HC^n$, where $\Ga < {\rm SU}(1,n)$ is a cocompact lattice isomorphic to the fundamental group of $M$. We then denote by $\map{\pi}{\HC^n}{M}$ the universal covering projection. Every compact complex-hyperbolic  manifold inherits a complex structure $I_{M}$ and a K\"ahler form $\om_M$ from $\HC^n$, given by $\pi^*\omega_M = \omega_{\HC^n}$. We define the normalized volume of $M$ by
\[
\Vol(M) \deq \int_{M} \om_{M}^{n}.
\]
 For $1$-forms $\alpha, \beta \in \Om^{1}(M;\hfr)$ taking values in some Lie algebra $\hfr$, we denote by $[\alpha \wedge \beta] \in \Om^{2}(M;\hfr)$ the Lie algebra valued $2$-form (often simply written $[\alpha, \beta]$ in the literature) defined by
\[
[\alpha \wedge \beta](v_1,v_2) \deq [\alpha(v_1), \beta(v_2)] - [\alpha(v_2), \beta(v_1)]
\]
for $v_{1},v_{2} \in TM$.

\begin{assumption}
	From now on we fix the following data: A compact complex-hyperbolic manifold $M = \Ga \backslash \HC^n$, a simple Hermitian Lie group $G$ with maximal compact subgroup $K$ of $G$ and corresponding Cartan decomposition $\gfr = \kfr \oplus \pfr$, an orientation and corresponding complex structure $I_G$ on the symmetric space $\X_G = G/K$,  a non-zero $G$-invariant $2$-form $\om_{G} \in \Om^{2}(\X_{G})^{G}$ on $\X_G$, and a representation $\map{\rho}{\Ga}{G}$. We further assume that $\rho$ is \emph{reductive}, i.e., that the Zariski closure of $\rho(\Gamma)$ is a reductive subgroup of $G$.
\end{assumption}

In the sequel we denote by $\map{\Ad}{G_\C}{{\rm Aut}(\gfr_\C)} < \GL(\gfr_\C)$ and  $\map{\ad}{\gfr_\C}{{\mathfrak {gl}}(\gfr_\C)}$ the adjoint representations of $G_\C$ and $\mathfrak g_\C$ respectively. If $H$ is a closed subgroup of $G_\C$ with Lie algebra $\mathfrak h$ and $\mathfrak{m} \subset \gfr_\C$ is an $H$-invariant subspace, then we denote the corresponding representations of $H$ and $\mathfrak h$ on $\mathfrak m$ by ${\Ad}_H^{\mathfrak m}$, respectively  ${\ad}_{\mathfrak h}^{\mathfrak m}$, in particular we denote
\begin{equation} \label{adkcpc}
	\map{\Ad_{K}^{\pfr_{\C}}}{K}{\GL(\pfr_{\C})} \quad \text{and} \quad \map{\ad_{\kfr}^{\pfr_{\C}}}{\kfr}{\GL(\pfr_{\C})}.
\end{equation}
Recall that the compact real form of $\gfr_{\C}$ is given by the Lie algebra $\ufr \deq \kfr \oplus i \, \pfr \subset \gfr_{\C}$, hence in particular $\gfr_{\C} = \ufr \oplus i \, \ufr$. We will denote by $\map{\tau}{\gfr_\C}{\gfr_\C}$ the complex conjugation with respect to $\ufr$ in $\gfr_{\C}$, i.e. $\tau(X + i\,Y) = X - i\,Y$ for $X,Y \in \ufr$.

\subsection{Higgs principal bundles}

\begin{definition}\label{GHiggsPB}
	By a \emph{$G$-Higgs principal bundle} over the manifold $M$ we mean a triple $(P,A,\vphi)$ consisting of
\begin{itemize}
	\item a principal $K$-bundle $P \to M$,
	\item a connection $1$-form $A \in \Om^{1}(P;\kfr)$,
	\item a $1$-form $\vphi \in \Omega^{1,0}(M, P \times_{\Ad} \pfr_\C)$ with values in the bundle $P \times_{\Ad} \pfr_\C$,
\end{itemize}
such that the pair $(A,\vphi)$ satisfies the \emph{complex Higgs equations}
\begin{equation} \label{ComplexHiggsEquations}
	\left\{
	\begin{aligned}
		\delbar_{A} \vphi = 0, \\
		[\vphi \wedge \vphi] = 0, \\
		F_{A} - [\vphi \wedge \tau(\vphi)] = 0.
	\end{aligned}
	\right.
\end{equation}
For an explanation of the notation see Remark \ref{RemarkComplexHiggsEquations2} below.
\end{definition}

The $1$-form $\vphi \in \Omega^{1,0}(M, P \times_{\Ad} \pfr_{\C})$ is often called the \emph{Higgs field}. If $M$ is a curve, the complex Higgs equations \eqref{ComplexHiggsEquations} are also known as \emph{Hitchin's equations} or \emph{self-duality equations} \cite{HitchinSelfDuality}. In this case, the second equation is automatically satisfied for dimensional reasons.

\begin{remark}
	Our definition of a $G$-Higgs principal bundle is adapted from the various definitions appearing in the literature, see e.g.~Hitchin \cite{HitchinSelfDuality}, Koziarz and Maubon \cite{KoziarzMaubon}, Sec.\,2 and Bradlow, Garc\'ia-Prada, and Gothen \cite{BGG}, Sec.\,2. Sometimes a $G$-Higgs principal bundle is alternatively defined as a holomorphic $K_{\C}$-principal bundle $P_{K_{\C}}$ over $M$, together with a holomorphic $1$-form $\th \in \Omega^{1,0}(M, P_{K_{\C}} \times_{\Ad} \pfr_\C)$ satisfying the equation $[\th \wedge \th] = 0$. However, as explained in \cite{KoziarzMaubon}, Sec.\,2 and \cite{BGG}, Sec.\,2.3, this definition is in fact equivalent to ours. The bundle $P_{K_\C}$ is simply given by the 
complexification $P_{K_{\C}} = P \times_{K} K_{\C}$ of the bundle $P_K$.
\end{remark}

\begin{remark} \label{RemarkComplexHiggsEquations2}
	Let us comment on the notation appearing in Definition \ref{GHiggsPB}.

(i) The twisted Dolbeault operator $\delbar_{A}$ is defined as the $(0,1)$-part of the connection
\[
\map{d_{A}}{\Omega^{1}(M, P \times_{\Ad} \pfr_\C)}{\Omega^{2}(M, P \times_{\Ad} \pfr_\C)}, \quad d_{A} \vphi \deq d \vphi + \left[ A \wedge \vphi \right]
\]
on the bundle $P \times_{\Ad} \pfr_{\C}$. More concretely, this means that
\[
\delbar_{A} \vphi = \frac{1}{2} \left( d_{A}\vphi + i \circ d_{A}\vphi \circ I_{M} \right).
\]
The first equation in \eqref{ComplexHiggsEquations} then says that the $1$-form $\vphi$ is holomorphic with respect to the holomorphic structure on the bundle $P \times_{\Ad} \pfr_{\C}$ defined by the connection $d_{A}$.

	(ii) The assumption that the $1$-form $\vphi \in \Omega^{1,0}(M, P \times_{\Ad} \pfr_{\C})$ be of type $(1,0)$ means the following. Denote by $I_{M}$ the complex structure on $M$ and by $i$ the standard complex structure on $\pfr_\C = \pfr \oplus i\,\pfr$. Then $\vphi$ is complex linear when viewed as a linear map $\map{\vphi}{TM}{P \times_{\Ad} \pfr_\C}$, i.e. $i \circ \vphi = \vphi \circ I_{M}$.

(iii) The curvature $2$-form $F_{A} \in \Omega^{2}(M, P \times_{\Ad} \kfr)$ is defined as
\[
F_{A} \deq dA + \frac{1}{2} \left[ A \wedge A \right].
\]

(iv) The conjugation $\tau: \gfr_\C \to \gfr_\C$ at the compact real form induces a fiberwise involution
\[
{P \times_{\Ad} \pfr_\C}\to {P \times_{\Ad} \pfr_\C}
\]
on the bundle $P \times_{\Ad} \pfr_\C$. By abuse of notation, we also denote this fiberwise involution by $\tau$.
\end{remark}
For more details on $G$-Higgs principal bundles and the complex Higgs equations, the reader may consult \cite{HitchinSelfDuality,KoziarzMaubon,BGG,Wentworth}. For background information on gauge theory we recommend \cite{Wehrheim}, App.\,A.

\subsection{Constructing Higgs bundles from representations} \label{SubSectionConstructingHiggs}

The non-abelian Hodge correspondence (see e.g.~\cite{HitchinSelfDuality,DonaldsonHarmonic,Corlette,SimpsonVariations,SimpsonLocalSystems,SimpsonModuliI,SimpsonModuliII,BGM,GarciaPradaGothenMundet/NAHC,Wentworth}) provides the link between reductive representations of the fundamental group of a given compact K\"ahler manifold into a semisimple Lie group on the one hand, and Higgs bundles over this K\"ahler manifold satisfying certain stability conditions on the other hand. For our purpose in this article it will be sufficient to spell out how a $G$-Higgs principal bundle $(P,A,\vphi)$ over $M$ is constructed from a reductive representation $\map{\rho}{\Ga}{G}$. In the following we will denote this assignment by
\[
\NAH_{G}\!: \rho \mapsto (P,A,\vphi)
\]
for short. The details are carefully explained in Koziarz and Maubon \cite{KoziarzMaubon}, Sec.\,2 (see also Bradlow, Garc\'ia-Prada, and Gothen \cite{BGG}, Sec.\,2). For later reference we briefly recall the main steps of this construction.

\bigskip

\noindent \textbf{Step 1: The principal $K$-bundle $P$}

\smallskip

Fix a reductive representation $\map{\rho}{\Ga}{G}$. Recall from Sect.\,\ref{SecToledoHiggsPrelims} that $\Ga < G$ acts on $\HC^n$ with quotient $M = \Ga \backslash \HC^n$, and that $\X_{G} = G/K$. We denote by
\[
P_{G} \deq \HC^n \times_{\rho} G \to M
\]
the flat principal $G$-bundle associated to the representation $\rho$. It is well-known (see e.g. \cite{Kobayashi}, Sec.~I.5) that there is a correspondence between the following three types of objects:
\begin{enumerate}[(i)]
\item principal $K$-subbundles $P$ of $P_{G}$;
\item metrics on the principal $G$-bundle $P_{G}$, i.e., reductions of its structure group to the maximal compact subgroup $K$ of $G$;
\item sections of the associated fiber bundle $P_{G} \times_{G} \X_{G} \cong \HC^n \times_{\rho} \X_{G} \to M$, or, equivalently, smooth $\rho$-equivariant maps $\map{f}{\HC^n}{\X_{G}}$.
\end{enumerate}
Here the $K$-subbundle $P$ of $P_{G}$ associated with a map $f$ as in (iii) is determined by the relation $\pi^{\ast}P = f^{\ast}G$, where $\map{\pi}{\HC^n}{M}$ denotes the projection and $G$ is thought of as the canonical principal $K$-bundle $G \to \X_G=G/K$. More concretely, the bundle $P$ is obtained as follows: since the map $f$ is $\rho$-equivariant, the principal $K$-bundle $f^{\ast}G \to \HC^{n}$ is invariant under the action of $\Ga$ from the left and hence descends to the principal $K$-bundle $P \to M$. Such a principal $K$-subbundle in turn defines a reduction of structure group for the bundle $P_G$. We say that a metric on $P_G$ is \emph{harmonic} if the corresponding map $f$ is harmonic. 

The non-abelian Hodge correspondence relies on the basic fact, due to Donaldson \cite{DonaldsonHarmonic} and Corlette \cite{Corlette} that, whenever the representation $\rho$ is reductive there exists a unique $\rho$-equivariant harmonic map $\map{f}{\HC^n}{\X_{G}}$ as in (iii) above. Accordingly, since our representation $\rho$ is assumed to be reductive we obtain a unique principal $K$-subbundle $P$ of $P_G$ corresponding to the harmonic map $f$. It fits into the following commutative diagram.
\[
\begin{tikzcd}
  & P_{G} \arrow[ddd] & & \HC^{n} \times G \arrow[ddd] \arrow[ll] \\
P \arrow[ddd] \arrow[ur, hook] & & f^{\ast}G \arrow[ll] \arrow[ur, hook] \arrow[ddd] \arrow[rr] & & G \arrow[ddd] \\
\\
  & M & & \HC^{n} \arrow[ll, near end, swap, "\pi"] \\
M \arrow[ur, "\mathrm{id}"] & & \HC^{n} \arrow[ur, "\mathrm{id}"] \arrow[ll] \arrow[rr, "f"] & & \X_{G} \\
\end{tikzcd}
\]

\bigskip

\noindent \textbf{Step 2: The connection form $A$}

\smallskip

Let $\alpha_G\in \Omega^1(G, \gfr)^{G}$ be the Maurer-Cartan form on $G$ and let $\map{p_{\kfr}}{\gfr}{\kfr}$ denote the canonical projection along $\pfr$ associated to the Cartan decomposition $\gfr = \kfr \oplus \pfr$. The $\kfr$-valued $1$-form $p_\kfr \circ \alpha_G \in \Omega^1(G, \kfr)$ is a connection $1$-form on the canonical principal $K$-bundle $G \to \X_{G}=G/K$. It pulls back to a connection $1$-form on the principal $K$-bundle $f^{\ast}G \to \HC^{n}$, and then further descends under the action of $\Ga$ to a connection $1$-form $A \in \Om^{1}(P,\kfr)$ on the principal $K$-bundle $P$.

\bigskip

\noindent \textbf{Step 3: The Higgs field $\vphi$}

\smallskip

Denote by $\map{p_{\pfr}}{\gfr}{\pfr}$ the canonical projection along $\kfr$ associated to the Cartan decomposition $\gfr = \kfr \oplus \pfr$. The $\pfr$-valued $1$-form $p_\pfr \circ \alpha_G \in \Omega^1(G, \pfr)$, where $\alpha_{G}$ is the Maurer-Cartan form as above, is horizontal and $K$-invariant, hence descends to a $1$-form $\alpha_G^\pfr \in \Omega^1(\X_G, G \times_{\Ad}\pfr)$ on the symmetric space $\X_G = G/K$ taking values in the bundle $G \times_{\Ad}\pfr$ of the canonical principal $K$-bundle $G \to \X_{G}=G/K$. Since $f$ is $\rho$-equivariant and $\alpha_G$ is $G$-invariant, the $1$-form $f^{\ast}\alpha_G^\pfr \in \Omega^1(\HC^{n}, f^{\ast}G \times_{{\Ad}}\pfr)$ is $\Gamma$-invariant and hence descends to a unique $1$-form $\psi\in \Omega^1(M, P \times_{{\Ad}}\pfr)$ on $M$ satisfying the relation
\begin{equation} \label{DefTheta}
	\pi^*\psi = f^*\alpha_G^\pfr.
\end{equation}
This $1$-form $\psi$ will be called the \emph{real Higgs field}. We then define the complex Higgs field $\vphi \in \Omega^{1,0}(M, P \times_{\Ad} \pfr_\C)$ as the $(1,0)$-part of $\psi$ with respect to the complex structure $I_{M}$ on $M$ and the standard complex structure $i$ on $\pfr_\C = \pfr \oplus i\,\pfr$. More concretely, this means that
\begin{equation} \label{EquationFormulavphi}
	\vphi \deq \frac{1}{2} \left( \psi - i \circ \psi \circ I_{M} \right).
\end{equation}
The next lemma relates the real and complex Higgs fields.

\begin{lemma} \label{LemmaRelationRealAndComplexHiggsField}
	The real Higgs field $\psi\in \Omega^1(M, P \times_{{\rm Ad}}\pfr)$ and the complex Higgs field $\vphi \in \Omega^{1,0}(M, P \times_{\Ad} \pfr_\C)$ are related by the formula
\[
\psi = \vphi - \tau(\vphi).
\]
\end{lemma}

\begin{proof}
	Recall from Remark \ref{RemarkComplexHiggsEquations2}\,(v) that $\map{\tau}{\gfr_\C}{\gfr_\C}$ denotes complex conjugation in $\gfr_{\C}$ with respect to $\ufr = \kfr \oplus i \, \pfr$. It follows that $\tau|_{\pfr} = -\Id$ and $\tau|_{i\pfr} = +\Id$. Since $\psi$ takes values in $\pfr$ we conclude that $\tau \circ \psi = -\psi$ and $\tau \circ i \circ \psi = i \circ \psi$. Thus we compute with \eqref{EquationFormulavphi} that
\[
\vphi - \tau(\vphi) = \frac{1}{2} \psi - \frac{1}{2} \, i \circ \psi \circ I_{M} - \frac{1}{2} \, \tau \circ \psi + \frac{1}{2} \, \tau \circ i \circ \psi \circ I_{M} = \psi.
\]
\end{proof}

\bigskip

\noindent \textbf{Step 4: The complex Higgs equations}

\smallskip

Harmonicity of the $\rho$-equivariant map $\map{f}{\HC^n}{\X_{G}}$ in Step 1 implies that the pair $(A,\vphi)$ does in fact satisfy the complex Higgs equations \eqref{ComplexHiggsEquations}. This was first observed by Donaldson \cite{DonaldsonHarmonic}. We refer the reader to Koziarz and Maubon \cite{KoziarzMaubon}, Sec.\,2.1 and 2.2 for further details.

\section{Toledo invariants in terms of Higgs principal bundles }

\subsection{A differential-geometric formula for the Toledo invariant}
Throughout this section we will freely use the notations introduced in Section \ref{SectionHodgeCorrespondence}. Moreover, we fix a non-zero invariant $2$-form $\om_{G} \in \Om^{2}(\X_{G})^{G}$ (subject to additional restrictions in Subsection \ref{SubSectionToledoForTraceForms}), and consider representations $\map{\rho}{\Ga}{G}$. Our starting point is the following expression, due to Burger and Iozzi, for $T_{\om_G}(\rho)$ in terms of purely differential-geometric data.

\begin{proposition}[\cite{BurgerIozzi/BoundedDifferentialForms}, Lemma 5.3] \label{PropositionBurgerIozziFormulaForToledoInvariant}
	Fix a representation $\map{\rho}{\Ga}{G}$ and a non-zero invariant $2$-form $\om_{G} \in \Om^{2}(\X_{G})^{G}$. Let $\map{f}{\HC^n}{\X_{G}}$ be a smooth $\rho$-equivariant map. Then the Toledo invariant of $\rho$ with respect to $\om_G$ is given by
\[
T_{\om_G}(\rho) = \frac{1}{\Vol(M)} \cdot \int_{M} (\pi_{\ast}f^{\ast}\om_G) \wedge \om_{M}^{n-1}.
\]
\end{proposition}

\begin{proof}
	Combining formula (1.2) and Lemma 5.3 in \cite{BurgerIozzi/BoundedDifferentialForms}, we may express the Toledo invariant as
\[
T_{\om_{G}}(\rho) = \frac{\big\langle \pi_{\ast}f^{\ast}\om_G, \om_{M} \big\rangle}{\langle \om_{M},\om_{M} \rangle},
\]
where $\langle\cdot,\cdot\rangle$ denotes the $L^{2}$-inner product on the space $\Om^{2}(M)$ of $2$-forms on $M$ given by
\[
\langle \alpha,\beta \rangle \deq \int_{M} \alpha \wedge \beta \wedge \om_{M}^{n-2}.
\]
It follows that
\[
\big\langle \pi_{\ast}f^{\ast}\om_G,\om_{M} - T_{\om_{G}}(\rho) \cdot \om_{M} \big\rangle = 0,
\]
which means that
\[
\int_{M} (\pi_{\ast}f^{\ast}\om_G) \wedge \om_{M}^{n-1} = \int_{M} T_{\om_G}(\rho) \cdot \om_{M}^{n} = T_{\om_G}(\rho) \cdot \Vol(M).
\]
\end{proof}

\subsection{The von-Wangen-zu-Geroldseck formula}
The restriction $(\om_{G})_{o}$ of our non-zero invariant $2$-form $\om_{G} \in \Om^{2}(\X_{G})^{G}$ to the tangent space $T_{o}\X_{G} \cong \pfr$ is a $K$-invariant real-valued antisymmetric bilinear form on $\pfr$. It thus gives rise to a fiberwise antisymmetric bilinear form $\Omtilde_{G}$ on the bundle $G\times_{\Ad}\pfr \to \X_G$ associated with the canonical principal $K$-bundle $G \to \X_{G} = G/K$. This bilinear form $\Omtilde_{G}$ is invariant under the action of $G$ from the left. Since $f$ is $\rho$-equivariant, it follows that the pullback $f^{\ast}\Omtilde_{G}$ to the bundle $f^{\ast}G \times_{\Ad}\pfr \to \HC^{n}$ is $\Gamma$-invariant, hence descends to a unique fiberwise antisymmetric bilinear form $\Om_{G}$ on the bundle $P \times_{\Ad}\pfr \to M$ satisfying
\begin{equation}\label{DefOmega}
	\pi^{\ast}\Om_{G} = f^{\ast}\Omtilde_{G}.
\end{equation}
Given a $1$-form $\beta \in \Om^{1}(M, P \times_\Ad \pfr)$ we may then define a real-valued $2$-form $\Om_{G}(\beta, \beta) \in \Om^2(M)$ by setting
\begin{equation*}
	\Om_{G}(\beta, \beta)(\xi_1, \xi_2) \deq \Om_{G}(\beta(\xi_1), \beta(\xi_2))
\end{equation*}
for $\xi_1, \xi_2 \in TM$.

\begin{proposition}[von-Wangen-zu-Geroldseck\footnote{Friedrich Ludwig Franz Reichsfreiherr von Wangen zu Geroldseck, 1727-1782, Prince-Bishop of Basel, Switzerland.} formula] \label{PropvWzG}
	Fix a reductive representation $\map{\rho}{\Ga}{G}$ and a non-zero invariant $2$-form $\om_{G} \in \Om^{2}(\X_{G})^{G}$. The Toledo invariant of $\rho$ with respect to $\om_G$ is given in terms of the associated $G$-Higgs principal bundle $\NAH_{G}(\rho) = (P,A,\vphi)$ over $M$ by %
\begin{equation} \label{vWzG}
	T_{\om_G}(\rho) = \frac{1}{{\rm Vol(M)}} \cdot \int_M \Om_{G}\left( \vphi - \tau(\vphi),  \vphi - \tau(\vphi) \right)  \wedge \om_{M}^{n-1}.
\end{equation}
\end{proposition}

\begin{proof}
	In view of Proposition \ref{PropositionBurgerIozziFormulaForToledoInvariant} and Lemma \ref{LemmaRelationRealAndComplexHiggsField} it suffices to show that
\begin{equation} \label{PushPullvWzG}
	\Om_{G}(\psi,\psi) = \pi_{\ast}f^{\ast}\om_G.
\end{equation}
Now by \eqref{DefTheta}, \eqref{DefOmega} and the definition of the bilinear form $\Om_{G}$ above we have
\[
\begin{aligned}
	\pi^{\ast} \big( \Om_{G}(\psi,\psi) \big) &= \pi^{\ast}\Om_{G}\,(\pi^{\ast}\psi,\pi^{\ast}\psi) \\
	&= f^{\ast}\Omtilde_{G} \big( f^{\ast}\alpha_{G}^{\pfr},f^{\ast}\alpha_{G}^{\pfr} \big) \\
	&= f^{\ast} \big( \Omtilde_{G}(\alpha_{G}^{\pfr},\alpha_{G}^{\pfr}) \big) \\
	&= f^{\ast} \big( (\om_G)_{o}(p_{\pfr} \circ \alpha_{G},p_{\pfr} \circ \alpha_{G}) \big),
\end{aligned}
\]
where the notation is as in Step 3 in Section \ref{SubSectionConstructingHiggs}. Thus in order to establish the relation \eqref{PushPullvWzG} it remains to show that
\[
\om_{G}  =  (\om_{G})_{o}(p_{\pfr} \circ \alpha_{G},p_{\pfr} \circ \alpha_{G}).
\]
However, both forms are $G$-invariant and the claimed identity certainly holds at the base point $o$ of $\X_{G}$.
\end{proof}

\subsection{The Toledo invariant for trace forms}
\label{SubSectionToledoForTraceForms}

So far our choice of the non-zero invariant $2$-form $\om_{G} \in \Om^{2}(\X_{G})^{G}$ was arbitrary. We will be able to obtain a more concrete formula for the Toledo invariant in terms of degrees of certain vector bundles once we restrict to a special class of such invariant $2$-forms.

For this let $V$ be a complex vector space and let $\map{\sigma}{K}{{\rm GL}(V)}$ be a representation with differential $\map{d\sigma}{\kfr}{\L{gl}(V)}$. We then define the associated \emph{trace form} as the alternating bilinear form on $\pfr$ given by
\begin{equation} \label{DefTraceFormo}
	(\om_{G}^{\sig})_o(X, Y) \deq \tr(d\sig(I) \circ d\sig([X,Y])) \quad (X,Y \in \pfr),
\end{equation}
where $I \in \kfr$ is the central element defining the complex structure $I_{G}$ on $\X_{G}$ via \eqref{DefI} and $\tr$ denotes the complex trace of a complex linear endomorphism.
\begin{lemma} \label{LemmaExtendTraceForm}
	For every complex representation $\map{\sigma}{K}{{\rm GL}(V)}$ with ${\rm tr}(d\sigma(I)^2) \neq 0$ the associated trace form $(\om_{G}^{\sig})_o$ is non-zero, takes real values and is invariant under $K$, hence extends to a non-zero invariant $2$-form $\om_{G}^{\sig} \in \Om^2(\X_{G})^{G}$.
\end{lemma}

\begin{proof} The key observation is that $d\sigma(I)$ is central in $d\sigma(\kfr)$ (cf. Section \ref{SubSectionSymmetricSpaces}). This, together with conjugation-invariance of the trace implies that $(\om_{G}^{\sig})_o$ is $K$-invariant. It also implies that the $d\sigma(I)$-eigenspaces $V_\pm$ are invariant under $d\sigma(\kfr)$. Hence, in order to show that the trace of $d\sigma(I) \circ d\sigma([X, Y])|_{V_\pm}$ is real it suffices to show that the trace of $d\sigma([X, Y])|_{V_\pm}$ is purely imaginary for all $X,Y \in \pfr$. However, this follows from the fact that since $K$ is compact, $\sigma(K)$ takes values in a maximal compact subgroup of ${\rm GL}(V)$, which is conjugate to $U(V)$. Thus $d\sigma(\kfr)$ takes values in a conjugate of $\mathfrak u(V)$, and we are left only with showing that  $(\om_{G}^{\sig})_o \neq 0$.

For this, we extend the representation $d\sigma$ complex-linearly to a representation $\map{d\sigma_\C}{\kfr_\C}{\L{gl}(V)}$. Then it suffices to show that the the complex-linear extension
\begin{equation} \label{CplxLinExtTraceForm}
	(\om_{G}^{\sig})_o^\C : \pfr_{\C} \times \pfr_{\C} \to \C, \quad (\om_{G}^{\sig})_o^\C (X, Y) \deq \tr(d\sig(I) \circ d\sig_\C([X,Y])) 
\end{equation}
of $(\om_{G}^{\sig})_o$ is non-zero. Recall that, by a classical theorem of Gordon Brown \cite{Brown} every element in a complex semisimple Lie algebra arises as a commutator. Thus there exists $X, Y \in \gfr_\C$ such that $[X, Y] = I$. Since ${\kfr}_\C = \mathfrak z(\kfr) \oplus [\kfr_\C, \kfr_\C]$ and $[\kfr_\C, \pfr_\C] \subset \pfr_\C$ we necessarily have $X,Y \in \pfr_\C$. For this choice of elements $X$ and $Y$ we then have
\[
\tr(d\sig(I) \circ d\sig([X,Y])) = \tr(d\sig(I)^2) \neq 0
\]
which shows that the bilinear form \eqref{CplxLinExtTraceForm} is non-zero.
\end{proof}

In the situation of Lemma \ref{LemmaExtendTraceForm} we refer to $\om_{G}^{\sig} \in \Om^2(\X_{G})^{G}$ as the \emph{invariant $2$-form associated with $\sigma$}.

\begin{assumption} \label{AssumptionPairVsigma}
For the rest of this subsection we fix a complex vector space $V$ and a representation $\map{\sigma}{K}{{\rm GL}(V)}$. We assume that there exist non-zero real number $\mu^{\pm} \in \R\setminus\{0\}$ and a splitting $V = V_{+} \oplus V_{-}$ such that $d\sig(I)$ acts on the summands $V_{\pm}$ by
\[
d\sig(I)|_{V_{\pm}} = \mu^{\pm} \cdot i \cdot \mathbbm{1}_{V_{\pm}}.
\]
\end{assumption}

As we shall see in Sections \ref{SectionToledoInvariantInTermsOfAdjointBundles} and \ref{SectionToledoInvariantInTermsOfAdmissible} below, this assumption will be satisfied whenever $\sigma$ belongs to one of the following two types of representations:
\begin{itemize}
	\item the adjoint representation $\map{\Ad_{K}^{\pfr_\C}}{K}{\GL(\pfr_{\C})}$ (see Section \ref{SubSectionAdjoint});
	\item the restriction to $K$ of an admissible representation $\map{\sig}{G}{\GL(V)}$ (see Section \ref{SecAdmRep}).
\end{itemize}
Note that every representation $\sigma$ which satisfies Assumption \ref{AssumptionPairVsigma} automatically satisfies the assumptions of Lemma \ref{LemmaExtendTraceForm}, hence the associated invariant $2$-form $\om_{G}^{\sig} \in \Om^2(\X_{G})^{G}$ is well-defined and non-zero. We are now going to express the Toledo invariant of this $2$-form in terms of degrees of certain vector bundles, which are constructed as follows.

Given a principal $K$-bundle $P \to M$, the representation $\map{\sigma}{K}{{\GL}(V)}$ gives rise to the \emph{associated vector bundle}
\[
\V_{\sig} \deq P \times_{\sig} V.
\]
It is a complex vector bundle over $M$. By Assumption \ref{AssumptionPairVsigma}, the $K$-module $V$ admits a splitting as a sum of $K$-modules $V = V_{+} \oplus V_{-}$, which in turn determines a splitting of the associated vector bundle
\[
\V_{\sig} = \V_{\sig}^{+} \oplus \V_{\sig}^{-},
\]
with the subbundles $\V_{\sig}^{+}$ and $\V_{\sig}^{-}$ defined by
\[
\V_{\sig}^{\pm} \deq P \times_{\sig} V_{\pm}.
\]

We denote the corresponding frame bundles by $\GL(\V_{\sig})$ and $\GL(\V_{\sig}^{\pm})$. A connection $1$-form $A \in \Om^{1}(P;\kfr)$ on the bundle $P$ then gives rise to connection $1$-forms on the associated frame bundles, denoted by
\[
A_{\sig} \in \Om^{1}(\GL(\V_{\sig});\gl(V)), \quad A_{\sig}^{\pm} \in \Om^{1}(\GL(\V_{\sig}^{\pm});\gl(V_{\pm})).
\]
The curvature $2$-form $F_{A_{\sig}}$ of $A_{\sig}$ (see Remark \ref{RemarkComplexHiggsEquations2}\,(iv)) may then be expressed in terms of the curvature $2$-forms of $A$ and $A_{\sig}^{\pm}$ by the relation
\begin{equation} \label{RelationCurvatures}
	F_{A_{\sig}} = d \sig(F_{A}) = \left( \begin{matrix} F_{A_{\sig}^{+}} & 0 \\ 0 & F_{A_{\sig}^{-}} \end{matrix} \right).
\end{equation}
Before we state our result we briefly review how Chern-Weil theory provides a gauge theoretic interpretation of the degree of the vector bundles $\V_{\sig}^{\pm}$. We define the \emph{degree} ${\bf deg}( \V_{\sig}^{\pm} )$ of the vector bundle $\V_{\sig}^{\pm}$ by
\begin{equation} \label{FormulaDegree}
	{\bf deg}( \V_{\sig}^{\pm} ) \deq = \int_{M} c_{1}( \V_{\sig}^{\pm} ) \wedge \om_{M}^{n-1},
\end{equation}
where $c_{1}$ denotes the first Chern class.  Note that for curves $M$ this definition agrees with the usual notion of degree for line bundles, and that in general ${\bf deg}( \V_{\sig}^{\pm} )  = {\bf deg}( {\bf det} \V_{\sig}^{\pm} )$ depends only on the associated determinant line bundle ${\bf det} \V_{\sig}^{\pm}$. Chern-Weil theory tells us that the first Chern class of $\V_{\sig}^{\pm}$ is given in terms of the curvature by
\begin{equation} \label{FormulaChernWeil}
	c_{1}( \V_{\sig}^{\pm} ) = \frac{i}{2\pi} \cdot F_{A_{\sig}^{\pm}}.
\end{equation}
Note that this holds independently of the choice of the connection $1$-form $A_{\sig}^{\pm}$. For more details on Chern-Weil theory the reader is referred to \cite{Dupont}.

\begin{theorem} \label{ToledoVectorBundles}
	Fix a reductive representation $\map{\rho}{\Ga}{G}$, and let $(V, \sigma)$ be as in Assumption \ref{AssumptionPairVsigma} with associated invariant $2$-form $\om_{G}^{\sig} \in \Om^{2}(\X_{G})^{G}$. Then the Toledo invariant of $\rho$ with respect to $\om_{G}^{\sig}$ is given in terms of the associated $G$-Higgs principal bundle $\NAH_{G}(\rho) = (P,A,\vphi)$ by the formula
\begin{equation*}
	T_{\om_{G}^{\sig}}(\rho) = - \frac{2\pi}{\Vol(M)} \cdot \left( \mu^{+} \cdot {\bf deg}(\V_{\sig}^{+}) + \mu^{-} \cdot {\bf deg}(\V_{\sig}^{-}) \right),
\end{equation*}
where $\V_{\sig}^{\pm} \deq P \times_{\sig} V_{\pm}$ are the associated vector bundles.
\end{theorem}

\begin{proof}
	Let $(P,A,\vphi) = \NAH_{G}(\rho)$, and let $\V_{\sig} \deq P \times_{\sig} V$ and $\V_{\sig}^{\pm} \deq P \times_{\sig} V_{\pm}$. The proof is in two steps.

\smallskip

\noindent \textbf{Step 1.} We claim that it follows from the von-Wangen-zu-Geroldseck Formula \eqref{vWzG} that the Toledo invariant is given by
\begin{equation} \label{vWzG+}
	T_{\om_{G}^{\sig}}(\rho) = -\frac{1}{{\rm Vol(M)}} \cdot \int_M  \tr \left( d\sig(I) \circ d\sig(F_A) \right) \wedge \om_{M}^{n-1},
\end{equation}
where $F_{A} \in \Omega^{2}(M, P \times_{\Ad} \kfr)$ denotes the curvature of $A$. Indeed, for $\xi_{1},\xi_{2} \in TM$ we have
\begin{eqnarray*}
	&& \Om_{G}\left( \vphi - \tau(\vphi),  \vphi - \tau(\vphi) \right) (\xi_{1},\xi_{2}) \\
	&=&  {\rm tr}\left( d\sigma(I) \circ d\sigma \left(\left[ (\vphi - \tau(\vphi))(\xi_{1}), (\vphi - \tau(\vphi))(\xi_{2}) \right]\right) \right) \\
	&=&  {\rm tr}\left( d\sigma(I) \circ d\sigma \left( \frac{1}{2} \left[ \left( \vphi - \tau(\vphi) \right) \wedge \left( \vphi - \tau(\vphi) \right) \right](\xi_{1},\xi_{2}) \right) \right).
\end{eqnarray*}
Since $[-\wedge -]$ is symmetric, $\tau$ is a Lie algebra homomorphism, and using the second of the Higgs equations \eqref{ComplexHiggsEquations} we have
\[
\left[ (\vphi - \tau(\vphi)) \wedge (\vphi - \tau(\vphi)) \right] = [\vphi \wedge \vphi] - 2 \, [\vphi \wedge \tau(\vphi)] + \tau \left( [\vphi \wedge \vphi] \right) = -2 \, [\vphi \wedge \tau(\vphi)].
\]
Applying the third of the Higgs equations \eqref{ComplexHiggsEquations} we therefore arrive at
\[
\Om_{G}\left( \vphi - \tau(\vphi),  \vphi - \tau(\vphi) \right) = -{\rm tr} \left( d\sigma(I) \circ d\sigma(F_A) \right).
\]
Formula \eqref{vWzG+} now follows from Formula \eqref{vWzG}.

\smallskip

\noindent \textbf{Step 2.} By Formula \eqref{RelationCurvatures} and Assumption \ref{AssumptionPairVsigma} we compute
\[
\begin{aligned}
	\tr \left( d\sig(I) \circ d\sig(F_A) \right) &= \tr \left( d\sig(I) \circ \left( \begin{matrix} F_{A_{\sig}^{+}} & 0 \\ 0 & F_{A_{\sig}^{-}} \end{matrix} \right) \right) \\
	&= i \cdot \tr \left( \begin{matrix} \mu^{+} \cdot F_{A_{\sig}^{+}} & 0 \\ 0 & \mu^{-} \cdot F_{A_{\sig}^{-}} \end{matrix} \right) \\
	&= i \cdot \left( \mu^{+} \cdot \tr(F_{A_{\sig}^{+}}) + \mu^{-} \cdot \tr(F_{A_{\sig}^{-}}) \right).
\end{aligned}
\]
Hence applying Chern-Weil theory and combining Formulas \eqref{vWzG+}, \eqref{FormulaDegree} and \eqref{FormulaChernWeil} we obtain
\[
\begin{aligned}
	T_{\om_{G}^{\sig}}(\rho) &= -\frac{i}{{\Vol(M)}} \cdot \int_{M} \left( \mu^{+} \cdot \tr(F_{A_{\sig}^{+}}) + \mu^{-} \cdot \tr(F_{A_{\sig}^{-}}) \right) \wedge \om_{M}^{n-1} \\
	&= -\frac{2\pi}{{\Vol(M)}} \cdot \int_{M} \left( \mu^{+} \cdot c_{1}(\V_{\sig}^{+}) + \mu^{-} \cdot  c_{1}(\V_{\sig}^{-}) \right) \wedge \om_{M}^{n-1} \\
	&= - \frac{2\pi}{\Vol(M)} \cdot \left( \mu^{+} \cdot \deg(\V_{\sig}^{+}) + \mu^{-} \cdot \deg(\V_{\sig}^{-}) \right),
\end{aligned}
\]
which proves the theorem.
\end{proof}

\section{Toledo invariants in terms of  the complexified isotropy bundle}
\label{SectionToledoInvariantInTermsOfAdjointBundles}

\subsection{Splitting of the complexified isotropy bundle}
\label{SubSectionAdjoint}
The goal of this section is to establish the universal Milnor--Wood type inequality in isotropy form. We will use all the notations introduced in the previous two sections. Recall that the isotropy representation of $K$ complexifies to a representation $ \map{\Ad_{K}^{\pfr_{\C}}}{K}{\GL(\pfr_{\C})}$, and that the complexified tangent space $\pfr_\C$ of $\X$ admits a $K$-invariant splitting
\begin{equation} \label{SplittingpC}
	\pfr_{\C} = \pfr_{+} \oplus \pfr_{-},
\end{equation}
where $\pfr_{\pm}$ denote the eigenspaces of $(I_\X)_o \otimes \C = \ad_{\kfr}^{\pfr_{\C}}(I)$ for the eigenvalues $\pm i$. We deduce that the pair $(V,\sig)$ defined by
\[
V = \pfr_{\C}, \quad \sig = \map{\Ad_{K}^{\pfr_{\C}}}{K}{\GL(\pfr_{\C})}
\]
satisfies Assumption \ref{AssumptionPairVsigma} with
\[
V_{\pm} = \pfr_{\pm}, \quad \mu^{\pm} = \pm 1.
\]

Given a principal $K$-bundle $P \to M$, the corresponding associated vector bundles are the adjoint bundles
\[
\V_{\rm iso} = P \times_{\Ad_{K}^{\pfr_{\C}}} \pfr_{\C}, \quad \V_{\rm iso}^{\pm} = P \times_{\Ad_{K}^{\pfr_{\C}}} \pfr_{\pm},
\]
which are related by
\begin{equation} \label{pfrpmsplitting}
	\V_{\rm iso} = \V_{\rm iso}^+ \oplus \V_{\rm iso}^-.
\end{equation}

\begin{lemma} \label{DegreeSplitting1}
	The determinant line bundle of the bundle $\V_{\rm iso}$ is trivial, and hence we have ${\bf deg}(\V_{\rm iso}) = 0$.
\end{lemma}

\begin{proof}
	We have to show that the structure group of $P \times_{{\Ad_{K}^{\pfr_{\C}}}} \pfr_{\C}$ can be reduced to ${\SL}_{n}(\pfr_{\C})$. For this it suffices to show that $\tr_{\pfr_{\C}}({\ad}_{\kfr}^{\pfr_{\C}}(X)) = 0$ for all $X \in \kfr$. By conjugation-invariance of the trace, it suffices to show this for all $T \in \L t$, where $\L t < \kfr$ is a maximal torus. Since $G$ is Hermitian we have $\rk_\R(G) =\rk_\R(K)$, i.e. $\L t_\C$ is a Cartan subalgebra of $\gfr_{\C}$. Let us call a root $\alpha$ of $\gfr_{\C}$ with respect to $\L t_\C$ a compact, respectively non-compact root if the corresponing root space is contained in $\kfr_{\C}$, respectively $\pfr_{\C}$. The key observation is that if $\alpha$ is a non-compact root, then so is $-\alpha$. This implies that $\tr_{\pfr_{\C}}({\ad}_{\kfr}^{\pfr_{\C}}(T)) = 0$ for all $T \in \L t_{\C}$, and consequently also $\tr_{\pfr_{\C}}({\ad}_{\kfr}^{\pfr_{\C}}(X)) = 0$ for all $X \in \kfr$.
\end{proof}

In terms of the splitting \eqref{pfrpmsplitting}, Lemma \ref{DegreeSplitting1} can be restated in the form
\begin{equation} \label{DegreeSplitting2}
	{\bf deg}( \V_{\rm iso}^+) + {\bf deg}( \V_{\rm iso}^-) = 0. 
\end{equation}

\subsection{Toledo invariant and Milnor--Wood type inequality}
\label{SubSectionCanoncialTraceForm}

By definition, the invariant $2$-form associated with the pair $(V,\sig) = (\pfr_{\C},\Ad_{K}^{\pfr_{\C}})$ is precisely the \emph{canonical $2$-form} $\om_{G}^{\rm can} \in \Om^{2}(\X_{G})^{G}$ from the introduction, which is uniquely determined by the formula
\[
(\om_{G}^{\mathrm{can}})_o(X,Y) = \tr \left( \ad_{\kfr}^{\pfr_{\C}}(I) \circ \ad_{\kfr}^{\pfr_{\C}}([X, Y]) \right) 
\]
for $X, Y \in \pfr$. Using Relation \eqref{DegreeSplitting2}, as an immediate consequence of Theorem~\ref{ToledoVectorBundles} we obtain the following formula for the Toledo invariant with respect to the canonical trace form.

\begin{corollary} \label{ToledoCanoncialTraceForm}
	Fix a reductive representation $\map{\rho}{\Ga}{G}$. Then the Toledo invariant of $\rho$ with respect to the canonical $2$-form $\om_{G}^{\rm can}$ is given in terms of the associated $G$-Higgs principal bundle $\NAH_{G}(\rho) = (P,A,\vphi)$ by the formula
\begin{equation} \label{vWzG1}
	T_{\om_{G}^{\rm can}}(\rho) = -\frac{4 \pi}{\Vol(M)} \cdot {\bf deg}(\V_{\rm iso}^+) = \frac{4 \pi}{\Vol(M)} \cdot {\bf deg}(\V_{\rm iso}^-).
\end{equation}
\end{corollary}

Combining Corollary \ref{ToledoCanoncialTraceForm} with Lemma \ref{LemmaMW} yields the following bounds of Milnor--Wood type on the degrees of the bundles $P \times_{\Ad} \pfr_{\pm}$.

\begin{corollary} \label{MilnorWoodAdjoint}
	Fix a reductive representation $\map{\rho}{\Ga}{G}$ with associated $G$-Higgs principal bundle $\NAH_{G}(\rho) = (P,A,\vphi)$. Then
\[
\bigl|{\bf deg}(\V_{\rm iso}^{\pm})  \bigr| \le \frac{{\rm Vol(M)}}{4\pi \cdot |K^{\rm min}_{\om_{G}^{\rm can}}|} \cdot \rk(G).
\]
\end{corollary}
Here $\om_{G}^{\rm can} \in \Om^{2}(\X_{G})^{G}$ is the canonical $2$-form, and the constant $K^{\rm min}_{\om_{G}^{\rm can}}$ may be computed by Formula \eqref{KMinFormula} using only Lie algebra data.

To illustrate this inequality we discuss the following example.

\begin{example}\label{FirstSUpq} Let $Q$ be a quadratic form of signature $(p,q)$. We may assume that $p < q$ and that $Q$ is of the form $Q := {\rm diag}(1, \dots, 1, -1, \dots, -1)$. Let $G := {\rm SU}(Q) = {\rm SU}(p,q)$ and denote by $\gfr$ its Lie algebra. According to \cite{BurgerIozziWienhardHermitian}, Example 3.1, we have
\[
p_\gfr = (p-1)2+(q-p)+2 = p+q.
\]
We have a Cartan decomposition $\gfr = \kfr\oplus \pfr$, where
\[
\kfr = \left\{\left(\begin{matrix}A & 0 \\ 0& D\end{matrix}\right)\mid A \in {\rm u}_p(\C),\; D \in  {\rm u}_q(\C),\; \tr(A) + \tr(D) = 0 \right\} \cong \L s (\L u(p) \oplus \L u(q))
\]
and
\[
\pfr =  \left\{\left(\begin{matrix}0 & B \\ C& 0\end{matrix}\right)\mid B \in M_{p,q}(\C),\; C \in M_{q,p}(\C), \; B = C^*  \right\} \cong M_{p,q}(\C).
\]
The element $I$ defining the complex structure is given by
\[
I = \frac i {p+q} \cdot \left(\begin{matrix}q\cdot \mathbbm{1}_p & 0 \\ 0& -p\cdot \mathbbm{1}_q\end{matrix}\right),
\]
We choose $X \in \pfr$ as
\[
X = \left(\begin{matrix}0 &B \\ B^\top& 0\end{matrix}\right), \; B =  \left(\begin{matrix} 1 & 0 &\dots & 0 & 0 &\dots& 0\\ 0 & 1 & & 0 & 0 &\dots& 0\\ \vdots&&\ddots &\vdots &\vdots&&\vdots \\ 0 & 0 &\dots& 1 & 0 &\dots& 0
\end{matrix}\right).
\]
Under ${\rm ad}(X)^2$, $\L p$ decomposes into three eigenspaces with eigenvalues $4$, $1$ and $0$ given by
\[
\pfr_4 = \left\{\left(\begin{matrix}0 & B_1 & 0\\ B_1^* & 0 & 0\\ 0& 0 & 0\end{matrix}\right) \mid B_1 = -B_1^* \in \L{u}(p)\right\},
\]
\[
\pfr_1 = \left\{\left(\begin{matrix}0 & 0 & B_2\\ 0 & 0 & 0\\ B_2^*& 0 & 0\end{matrix}\right)\mid B_2 \in M_{q-p, p}(\C)\right\},
\]
\[
\pfr_0 =\left\{\left(\begin{matrix}0 & B_1 & 0\\ B_1^* & 0 & 0\\ 0& 0 & 0\end{matrix}\right)\mid B_1 = B_1^* \in M_{p, p}(\C)\right\}.
\]
Similarly, $\L k$ decomposes into three eigenspaces with eigenvalues $4$, $1$ and $0$ given by
\[
\kfr_4 = \left\{\left(\begin{matrix}A & 0 & 0\\ 0 & -A & 0\\ 0& 0 & 0\end{matrix}\right)\mid A \in \L u_p(\C)\right\},
\]
\[
\kfr_1 = \left\{\left(\begin{matrix}0 & 0 & 0\\ 0 & 0 & D_{12}\\ 0& -D_{12}^* & 0\end{matrix}\right)\mid D_{12} \in M_{q-p, p}(\C) \right\},
\]
\[
\kfr_0 = \left\{\left(\begin{matrix}A & 0 & 0\\ 0 & A & 0\\ 0& 0 & D_{22}\end{matrix}\right)\mid A \in \L u_p(\C), D_{22} \in M_{q-p. q-p}(\C), 2\tr(A) + \tr{D_{22}} = 0\right\}
\]
We thus get
\[
\tr({\rm ad}(X)^2) = 4\cdot p^2+1\cdot 2(q-p)p+4\cdot p^2+1\cdot 2(q-p)p = 4p^2+4pq.
\]
On the other hand, 
\[
[X, [I, X]] = \left[\left(\begin{matrix}0 &B \\ B& 0\end{matrix}\right),\left(\begin{matrix}0 &iB^{\top} \\ -iB^{\top}& 0\end{matrix}\right)\right] = \left( \begin{matrix} -2i \mathbbm{1}_p&&\\&2i \mathbbm{1}_p&\\ &&0
\end{matrix}\right),
\]
and thus
\[
(\ad_\kfr^\pfr([X,[I,Y]] \circ \ad_\kfr^\pfr(I))\left(\left(\begin{matrix} 0 & B_1 & B_2 \\ C_1 & 0 & 0 \\ C_2 & 0 & 0\end{matrix}\right)\right) = \left(\begin{matrix} 0 & 4B_1 & 2B_2 \\ 4C_1 & 0 & 0 \\ 2C_2 & 0 & 0\end{matrix}\right).
\]
Hence the real trace of $\ad_\kfr^\pfr([X,[I,Y]] \circ \ad_\kfr^\pfr(I)$ is given by
\[
\tr_\pfr(\ad_\kfr^\pfr([X,[I,Y]] \circ \ad_\kfr^\pfr(I)) = 4\cdot 2p^2 + 2\cdot 2 p(q-p) = 4p^2+4pq,
\]
and this coincides with the complex trace of $\ad_\kfr^{\pfr_\C}([X,[I,Y]] \circ \ad_\kfr^{\pfr_\C}(I)$. Plugging this into  \eqref{KMinFormula} we obtain
\begin{equation*} \label{KMinCase1}
	K^{\rm min}_{\om_{G}^{\rm can}} =  - \frac 1{p+q} \cdot \frac{4p^2+4pq} {{4p^2+4pq}} = -\frac{1}{p+q}.
\end{equation*}
The Milnor--Wood type inequality for the adjoint bundles from Corollary \ref{MilnorWoodAdjoint} in the case $G = {\rm SU}(p,q)$ therefore becomes
\[
\bigl|{\bf deg}(\V_{\rm iso}^{\pm})  \bigr| \le \frac{\Vol(M)}{4\pi } \cdot (p+q) \cdot \min\{p,q\}.
\]
\end{example}

\section{Toledo invariants in terms of vector bundles associated with admissible representations}
\label{SectionToledoInvariantInTermsOfAdmissible}

\subsection{Admissible representations} \label{SecAdmRep}

In \cite{BGG}, Bradlow, Garc{\'{\i}}a-Prada, and Gothen express the Toledo invariant of surface group representations in terms of the degrees of certain vector bundles. For example, they use the standard representation $\SU(p,q) \to \GL_{p+q}(\C)$ to associate with every $\SU(p,q)$-Higgs principal bundle a vector bundle. This vector bundle then naturally splits into two subbundles, and the authors express the Toledo invariant in terms of the degrees of these subbundles. In this section, we will provide a framework, based on the notion of \emph{admissible representation}, that generalizes this example.

Let us fix a complex vector space $V$, and consider a non-trivial representation $\map{\sig}{G}{\GL(V)}$.
We denote by
\begin{equation} \label{Mapdsigma}
	\map{d\sigma}{\gfr_\C}{\gl(V)}
\end{equation}
the complex linear extension of the linearization $\map{d\sigma}{\gfr}{\gl(V)}$ to the complexification $\gfr_\C$. While in general the complexification of a simple Lie algebra is merely semisimple, the complexified Lie algebra $\gfr_\C$ happens to be simple (as follows from the classification in \cite{Helgasson}). Since the representation \eqref{Mapdsigma} is non-trivial, it is therefore automatically faithful. Moreover, since $\gfr_{\C}$ is perfect, i.e. $[\gfr_{\C},\gfr_{\C}] = \gfr_{\C}$, the representation \eqref{Mapdsigma} takes values in $[\gl(V),\gl(V)] = \sl(V)$.

We claim that the endomorphism $d\sigma(I)$ of $V$ is diagonalizable. To see this, we first observe that $\map{\ad(I)}{\gfr_{\C}}{\gfr_{\C}}$, $X \mapsto [I,X]$ is diagonalizable. In fact, since $I$ is central in $\kfr$ we have $\ad(I)|_{\kfr} = 0$; moreover, $\ad(I)|_{\pfr_{\pm}} = \pm i \Id$ for the splitting $\pfr_\C = \pfr_+ \oplus \pfr_-$ as in \eqref{SplittingpC}. It then follows from Cor.\,6.4 in \cite{Humphreys} that $d\sigma(I)$ is diagonalizable. If $V_\lambda$ denotes the eigenspace of $d\sigma(I)$ of eigenvalue $\lambda$, then $d\sigma(\kfr_\C)$ preserves $V_\lambda$ (since $I$ is central in $\kfr_\C$) and 
\begin{equation} \label{WeightRoot}
	d\sigma(\pfr_{\pm}).V_{\lambda} \subset V_{\lambda\pm i}.
\end{equation}
Since ${\rm tr}(d\sigma(I)) = 0$, the weighted sum of the eigenvalues of $d\sigma(I)$ equals $0$, and in particular $d\sigma(I)$ has at least two eigenvalues since $d\sigma$ is faithful. Note that, in principle, the number of eigenvalues of $d\sigma(I)$ can be arbitrarily large.

\begin{definition}\label{DefAdmissibleRep}
	The representation $\map{\sig}{G}{\GL(V)}$ is called \emph{admissible} if the endomorphism $d\sigma(I) \in \gl(V)$ has exactly two eigenvalues.
\end{definition}

The following proposition shows that admissible representations share many structural properties with the standard representation of ${\rm SU}(p,q)$, and in particular, satisfy Assumption \ref{AssumptionPairVsigma}.

\begin{proposition} \label{PropAdmissability}
	Let $\map{\sig}{G}{\GL(V)}$ be an admissible representation. Then the following hold.
\begin{enumerate}
	\item[(i)] $\sig$ takes values in $\SL(V)$.
	\item[(ii)] There exist a splitting $V = V_{+} \oplus V_{-}$ such that
\[
d\sig(\kfr_{\C}) \subset \left( \begin{matrix} * & 0 \\ 0 & *\end{matrix} \right), \quad
d\sig(\pfr_{+}) \subset \left( \begin{matrix}  0 & * \\ 0 & 0 \end{matrix} \right), \quad
d\sig(\pfr_{-}) \subset \left( \begin{matrix}  0 & 0 \\ * & 0 \end{matrix} \right)
\]
and
\begin{equation} \label{dsigmaI}
	d\sig(I) = \frac{i}{\dim V} \cdot \left( \begin{matrix} \dim V_{-} \cdot \mathbbm{1}_{V_{+}} & 0 \\ 0 & - \dim V_{+} \cdot \mathbbm{1}_{V_{-}} \end{matrix} \right),
\end{equation}
\end{enumerate}
In particular, the pair $(V,\sig)$ satisfies Assumption \ref{AssumptionPairVsigma}.
\end{proposition}

\begin{proof}
	Assertion (i) follows from the observation above that the representation \eqref{Mapdsigma} takes values in $[\gl(V),\gl(V)] = \sl(V)$. To see (ii), we denote by $\lambda_{\pm}$ the two eigenvalues of $d\sig(I)$ and by $V_{\pm}$ the corresponding eigenspaces. Since $I$ is central in $\kfr_\C$, the latter preserves the splitting $V = V_{+} \oplus V_{-}$. Since $d\sig$ is faithful, $
\pfr_{\pm}$ act non-trivially on $V$. It thus follows from \eqref{WeightRoot} that we can ensure (by exchanging $\lambda_{\pm}$ if necessary) that 
\begin{equation} \label{EigenvalueDifferenceI}
	\lambda_{+} - \lambda_{-} = i.
\end{equation}
Then $d\sig(\pfr_{\pm}) V_{\pm} = \{0\}$ and $d\sig(\pfr_{\pm}) V_{\mp} = V_{\pm}$. This proves that $d\sig(\kfr_{\C})$ and $d\sig(\pfr_{\pm})$ are of the desired form. Moreover, since $d\sig$ takes values in $\mathfrak{sl}(V)$ we have
\[
\dim(V_{+}) \cdot \lambda_{+} + \dim(V_{-}) \cdot \lambda_{-} = \tr (d\sig(I)) = 0.
\]
Combining this with \eqref{EigenvalueDifferenceI} and using $\dim(V_{+}) + \dim(V_{-}) = \dim(V)$ we obtain
\[
\lambda_{\pm} =  \pm \frac{\dim V_{\mp}}{\dim V} \cdot i,
\]
which shows that $d\sig(I)$ is of the desired form.
\end{proof}

\subsection{Toledo invariant and Milnor--Wood type inequality}

Fix an admissible representation $\map{\sig}{G}{\GL(V)}$. By Proposition \ref{PropAdmissability} the pair $(V = V_{+} \oplus V_{-},\sig)$ satisfies Assumption \ref{AssumptionPairVsigma} with
\[
\mu^{\pm} = \pm \frac{\dim V_{\mp}}{\dim V}.
\]
Given a principal $K$-bundle $P \to M$, keeping the notation from Section \ref{SubSectionToledoForTraceForms} the corresponding associated vector bundles will be denoted by
\[
\V_{\sig} = P \times_{\sig} V, \quad \V_{\sig}^{\pm} = P \times_{\sig} V_{\pm}.
\]
As a consequence of Theorem~\ref{ToledoVectorBundles} we obtain the following formula for the Toledo invariant with respect to the invariant $2$-form $\om_{G}^{\sig}$ associated with the pair $(V, \sig)$.

\begin{corollary} \label{ToledoAdmissibleTraceForm}
	Fix a reductive representation $\map{\rho}{\Ga}{G}$ and an admissible representation $\map{\sig}{G}{\GL(V)}$. Then the Toledo invariant of $\rho$ with respect to the invariant $2$-form $\om_{G}^{\sig}$ is given in terms of the associated $G$-Higgs principal bundle $\NAH_{G}(\rho) = (P,A,\vphi)$ by the formula
\begin{equation} \label{vWzG2}
	T_{\om_{G}^{\sig}}(\rho) = -\frac{2 \pi}{\Vol(M)} \cdot {\bf deg}(\V_{\sig}^{+}) = \frac{2 \pi}{\Vol(M)} \cdot {\bf deg}(\V_{\sig}^{-}),
\end{equation}
where $\V_{\sig}^{\pm} \deq P \times_{\sig} V_{\pm}$ are the associated vector bundles.
\end{corollary}

\begin{proof}
	Since $\sig$ takes values in $\SL(V)$ by Proposition \ref{PropAdmissability}\,(i), it follows that ${\bf deg}(\V_{\sig}) = 0$. In view of the splitting $\V_{\sig} = \V_{\sig}^{+} \oplus \V_{\sig}^{-}$ this implies that
\[
{\bf deg}(\V_{\sig}^{+}) + {\bf deg}(\V_{\sig}^{-}) = 0.
\]
Using this relation, we obtain from Theorem~\ref{ToledoVectorBundles} that
\[
\begin{aligned}
	T_{\om_{G}^{\sig}}(\rho) &= - \frac{2\pi}{\Vol(M)} \cdot \left( \mu^{+} \cdot {\bf deg}(\V_{\sig}^{+}) + \mu^{-} \cdot {\bf deg}(\V_{\sig}^{-}) \right) \\
	&= - \frac{2\pi}{\Vol(M)} \cdot \left( \frac{\dim V_{-}}{\dim V} \cdot {\bf deg}(\V_{\sig}^{+}) + \frac{\dim V_{+}}{\dim V} \cdot {\bf deg}(\V_{\sig}^{+}) \right) \\
	&= -\frac{2 \pi}{\Vol(M)} \cdot {\bf deg}(\V_{\sig}^{+}),
\end{aligned}
\]
and likewise for the second identity.
\end{proof}

Combining Corollary \ref{ToledoAdmissibleTraceForm} with Lemma \ref{LemmaMW} yields the following bounds of Milnor--Wood type on the degrees of the vector bundles $\V_{\sig}^{\pm} = P \times_{\sig} V_{\pm}$.

\begin{corollary} \label{MilnorWoodAdmissible}
	Fix a reductive representation $\map{\rho}{\Ga}{G}$ with associated $G$-Higgs principal bundle $\NAH_{G}(\rho) = (P,A,\vphi)$, and an admissible representation $\map{\sig}{G}{\GL(V)}$. Then
\[
\bigl| \deg(\V_{\sig}^{\pm}) \bigr| \le \frac{\Vol(M)}{2\pi \cdot |K^{\rm min}_{\om_{G}^{\sig}}|} \cdot \rk(G).
\]
\end{corollary}
Here $\om_{G}^{\sig} \in \Om^{2}(\X_{G})^{G}$ is the invariant $2$-form associated with the pair $(V, \sig)$, and the constant $K^{\rm min}_{\om_{G}^{\sig}}$ may be computed by Formula \eqref{KMinFormula} using only Lie algebra data. We illustrate the corollary by working out the case of the standard representation of ${\rm SU}(p,q)$.

\begin{example}
As in Example \ref{FirstSUpq} we consider the group ${\rm SU}(p,q)$, $p<q$, and we use the same notations as introduced there. Moreover, we denote by
 $\sigma: {\rm SU}(p,q) \to {\rm GL}_{p+q}(\C)$ the standard representation. Recall from Example \ref{FirstSUpq} that $p_\gfr = p+q$, and that we can choose $X \in \pfr$ such that
\[
\tr({\rm ad}(X)^2) = 4p^2+4pq
\]
and
\[
[X, [I, X]] = \left( \begin{matrix} -2i\mathbbm{1}_p&&\\&2i\mathbbm{1}_p&\\ &&0
\end{matrix}\right).
\]
We deduce that
\begin{eqnarray*}
\omega^\sigma_o(X, [I,X]) = {\rm tr}\left(\frac{i}{p+q}\left( \begin{matrix} q \cdot\mathbbm{1}_p & & \\  & p  \cdot \mathbbm{1}_p& \\ &&p \cdot \mathbbm{1}_{q-p} \end{matrix}\right)\left(\begin{matrix} -2i\mathbbm{1}_p&&\\&2i\mathbbm{1}_p&\\ &&0
\end{matrix}\right)\right) = 2p,
\end{eqnarray*}
and hence, by \eqref{KMinFormula}, 
\[
|K^{\rm min}_{\om_{G}^{\sig}}| = \frac{1}{p_\gfr} \cdot \frac{\tr({\rm ad}(X)^2)}{(\omega_G)_o(X,[I,X])} = \frac{1}{p+q} \cdot \frac{4p^2+4pq}{2p} = 2.
\]
The Milnor--Wood type inequality from Corollary \ref{MilnorWoodAdmissible} for the vector bundles associated with the standard representation of $G = \SU(p,q)$ therefore becomes
\begin{equation}
	\bigl|{\bf deg}({\bf V}_{\sig}^{\pm})  \bigr| \le \frac{{\rm Vol(M)}}{4\pi} \cdot {\rm min}\{p,q\}.
\end{equation}
In the special case $\dim(M)=1$, in which $M = \Sigma_{g}$ is a Riemann surface of genus $g \geq 2$ this inequality specializes to
\[
\bigl| {\bf deg}(\V_{\sig}^{\pm}) \bigr| \le \frac{|\chi(\Sigma_{g})|}{2} \cdot {\rm min}\{p,q\} = (g-1) \cdot {\rm min}\{p,q\},
\]
which agrees with the inequality obtained in \cite{BGG}.
\end{example}

\end{document}